\documentclass[11pt,a4paper]{article}

\usepackage[colorlinks]{hyperref}

\usepackage[table]{xcolor}
\usepackage{slashed}
\usepackage{xkeyval}
\usepackage{url,amsmath,ifthen,amssymb,latexsym,pstricks,mathrsfs,comment,amsthm,graphicx,tikz,tikz-cd,enumerate,accents,pgffor,cite,wrapfig,multicol,float}

\usepackage{bibspacing}

\usepackage{enumitem}
\usepackage[margin=10pt,font=small,labelfont=bf, labelsep=period]{caption}

\usepackage[T1]{fontenc}
\usepackage{wasysym}
\usepackage{stmaryrd}

\usepackage{geometry} 

\begin{document}

\newcommand{\nc}{\newcommand}
\nc{\rnc}{\renewcommand}

\nc\la\langle
\nc\ra\rangle
\nc\pres[2]{\la#1:#2\ra}
\nc\es\varnothing
\nc\bone{{\bf 1}}
\nc\si\sigma

\nc\set[2]{\{#1:#2\}}

\nc\ben{\begin{enumerate}[label=\textup{(\roman*)},leftmargin=7mm]}
\nc\bena{\begin{enumerate}[label=\textup{(\alph*)},leftmargin=7mm]}
\nc\een{\end{enumerate}}

\nc{\N}{\mathbb N}
\rnc{\P}{\mathbb P}
\nc\A{\mathcal A}
\nc\im{\operatorname{im}}

\nc\ba{{\bf a}}
\nc\bb{{\bf b}}
\nc\bc{{\bf c}}

\let\oldproofname=\proofname
\rnc{\proofname}{\rm\bf{\oldproofname}}

\nc\mt\mapsto
\nc\sub\subseteq

\nc{\COMMA}{,\qquad}
\nc{\COMMa}{,\ \ \ }
\nc{\OR}{\quad\text{or}\quad}
\nc{\AND}{\qquad\text{and}\qquad}
\nc{\ANd}{\quad\text{and}\quad}

\nc{\bit}{\begin{itemize}}
\nc{\eit}{\end{itemize}}
\nc{\bmc}{\begin{multicols}}
\nc{\emc}{\end{multicols}}
\nc{\pfitem}[1]{\medskip \noindent #1}
\nc{\firstpfitem}[1]{#1}

\numberwithin{equation}{section}

\newtheorem{thm}[equation]{Theorem}
\newtheorem{qu}[equation]{Question}
\newtheorem{lemma}[equation]{Lemma}
\newtheorem{cor}[equation]{Corollary}
\newtheorem{prop}[equation]{Proposition}
\newtheorem{conj}[equation]{Conjecture}
\newtheorem{hope}[equation]{Hope}

\theoremstyle{definition}

\newtheorem{defn}[equation]{Definition}
\newtheorem{rem}[equation]{Remark}
\newtheorem{eg}[equation]{Example}
\newtheorem{prob}[equation]{Problem}

\nc\pf{\begin{proof}}
\nc\epf{\end{proof}}

\title{Presentations for $\P^K$}

\date{}
\author{James East\footnote{Supported by ARC Future Fellowship FT190100632.}\\~\\
{\small Centre for Research in Mathematics and Data Science,} \\
{\small Western Sydney University, Locked Bag 1797, Penrith NSW 2751, Australia.} \\ 
{\small {\tt j.east\,@\,westernsydney.edu.au}}}

\maketitle

\begin{abstract}
It is a classical result that the direct product $A\times B$ of two groups is finitely generated (finitely presented) if and only if $A$ and $B$ are both finitely generated (finitely presented).  This is also true for direct products of monoids, but not for semigroups.  The typical (counter)example is when $A$ and $B$ are both the additive semigroup $\P=\{1,2,3,\ldots\}$ of positive integers.  Here $\P$ is freely generated by a single element, but $\P^2$ is not finitely generated, and hence not finitely presented.  In this note we give an explicit presentation for~$\P^2$ in terms of the unique minimal generating set; in fact, we do this more generally for~$\P^K$, the direct product of arbitrarily many copies of $\P$.

\emph{Keywords}: Semigroups, presentations, direct products, integers.

MSC2020: 20M05, 20M14.

\end{abstract}

\section{Introduction}\label{sect:intro}

Presentations (by generators and relations) are essential tools in algebra, and in many other parts of mathematics and science.  See for example \cite{LZ2015,Artin1947,Jones1983_2}; the introduction to \cite{East2020} contains many more references, and a more detailed discussion.  Similarly, methods for constructing new structures from existing ones abound, including (semi)direct products, wreath products and various extensions.  It is therefore of considerable interest to know how presentations for such constructions behave.  The current note concerns direct products of semigroups, and the intriguing complication caused by the lack of an identity element, as noted for example in \cite{HR1994,RRW1998}.

If groups $A$ and $B$ have presentations $\pres{X_1}{R_1}$ and $\pres{X_2}{R_2}$, then the direct product $A\times B$ has presentation $\pres{X_1\cup X_2}{R_1\cup R_2\cup R_3}$, where $R_3$ consists of all relations of the form $xy=yx$ with $x\in X_1$ and $y\in X_2$.  See for example \cite{J1980,MKS1966}.  Together with the fact that $A$ and $B$ are both homomorphic images of $A\times B$, it follows that $A\times B$ is finitely presented if and only if $A$ and $B$ both are.

The previous paragraph is true for direct products of monoids as well \cite{L1998,HR1994}, but not in general for semigroups.  The example often cited for this fact is when $A$ and~$B$ are both taken to be the additive semigroup of positive integers, $\P=\{1,2,3,\ldots\}$.  Indeed,~$\P$ is freely generated by $1$, so has finite (semigroup) presentation $\pres x\es$, yet $\P^2=\P\times\P$ is not even finitely generated.  To see this, note for example that any element of $\P^2$ of the form $(a,1)$ cannot be written as a sum of two elements of $\P^2$, and hence must be included in any generating set.

Thus, presentations for direct products of semigroups are far less `well-behaved' as for groups or monoids, and can even be somewhat `wild'.  Nevertheless, a number of interesting results are known \cite{RRW1998,MR2018,MR2019,CR2020}.  For example, in \cite{RRW1998}, necessary and sufficient conditions are given for a direct product $A\times B$ of semigroups to be finitely generated/presented.  These conditions of course include finite generation/presentability of $A$ and $B$, but more is required, including notions such as stability (whose definition we will not recall here).  Presentations for direct products are also constructed in \cite{RRW1998} when these conditions are satisfied, but not in general.  In \cite{HR1994}, the general situation was discussed with some pessimism, and the view was expressed that presentations for $A\times B$ (in general) were likely to be `at least as complicated as the multiplication table'.

In light of this, it seems to be of considerable interest to investigate specific instances of non-finitely-presented direct products of finitely-presented semigroups.  Moreover, it seems especially worthwhile to look at the case of $\P^2$, because of the clashing intuitions involved:  on the one hand,~$\P^2$ of course seems such a `tame' semigroup; but on the other hand, it arguably has the maximum possible degree of `wildness', as discussed above.

The purpose of the current note, therefore, is to give an explicit presentation for $\P^2$.  In fact, we prove a much more general result: Theorem \ref{t1} below gives a presentation for~$\P^K$, the direct product of arbitrarily many copies of $\P$.  The special case of $\P^2$ is discussed in Remark~\ref{r2}.  

Before we begin, we briefly note one additional strand of motivation for the current work.  When~$K$ is infinite, the semigroup $\P^K$ is of course uncountable.  Presentations for uncountable semigroups (or groups or monoids) are rarely discussed, as any generating set for such a semigroup~$S$ must have the same cardinality as $S$ itself.  However, if there is a `canonical' generating set (or family of such sets), it does seem of value to understand presentations with respect to such a set.  This is certainly the case for $\P^K$, as shown in Proposition \ref{p1} below.

We now take the opportunity to fix the notation we will be using for presentations.  For more background on semigroups, the reader is referred to a monograph such as \cite{CPbook} or \cite{Howie}.
The \emph{free semigroup} over a set $X$, denoted $X^+$, consists of all non-empty words over $X$ under concatenation.  Given a set $R\sub X^+\times X^+$ of pairs of words, we denote by $R^\sharp$ the congruence on~$X^+$ generated by~$R$.  A semigroup $S$ has \emph{presentation} $\pres XR$ if $S\cong X^+/R^\sharp$, or equivalently if there is a surmorphism (surjective homomorphism) $X^+\to S$ with kernel $R^\sharp$.  If $\phi$ is such a surmorphism, we say $S$ has \emph{presentation $\pres XR$ via $\phi$}.  The elements of $X$ and $R$ are referred to as \emph{generators} and \emph{relations}, respectively, and a relation $(u,v)\in R$ is typically displayed as an equation: $u=v$.

\section{The presentation}

We write $\P=\{1,2,3,\ldots\}$ and $\N=\{0,1,2,3,\ldots\}$ for the sets of positive and non-negative integers, respectively.  So $\P$ is an additive semigroup, and $\N$ its monoid completion.  We also fix an arbitrary set~$K$, and to avoid trivialities we assume that $|K|\geq2$.
Formally, the \emph{direct product}~$\P^K$ consists of all functions $K\to\P$, under point-wise addition.  Such a function will be identified with a $K$-tuple in the usual way.  

We adopt the convention that the entries of a tuple from $\P^K$ are denoted by the same letter as the tuple, so $\ba=(a_k)_{k\in K}$, $\bb=(b_k)_{k\in K}$, and so on.  We also write $\bone=(1)_{k\in K}$ for the $K$-tuple consisting entirely of $1$s.  Note that $\P^K=\bone+\N^K$, so that $\P^K$ is a principle ideal of the monoid~$\N^K$.

We call $\ba\in\P^K$ an \emph{atom} if there do not exist $\bb,\bc\in\P^K$ such that $\ba=\bb+\bc$.  

\begin{lemma}\label{l1}
An element $\ba\in\P^K$ is an atom if and only if $a_k=1$ for some $k\in K$.
\end{lemma}

\pf
If $a_k\geq2$ for all $k$, then $\ba-\bone\in\P^K$, and from $\ba=(\ba-\bone)+\bone$ it follows that $\ba$ is not an atom.

Conversely, if $\ba=\bb+\bc$ for some $\bb,\bc\in\P^K$, then $a_k=b_k+c_k\geq2$ for all $k$.
\epf

Atoms will play a crucial role in all that follows.  We will write $\A$ for the set of all atoms of~$\P^K$.

Because $\P$ is well ordered, we have a function $\mu:\P^K\to\P$, defined for $\ba\in\P^K$ by
\[
\mu(\ba)=\min_{k\in K}a_k.
\]

\begin{lemma}\label{l2}
If $\ba\in\P^K$, then $\ba=m\bone+\bb$ for unique $m\in\N$ and $\bb\in \A$.
\end{lemma}

\pf
To establish existence, put $m=\mu(\ba)-1$ and $\bb=\ba-m\bone$.  Certainly $\ba=m\bone+\bb$, and since $\mu(\ba)\in\P$ we have $m\in\N$.  By definition of $\mu(\ba)$, we have $a_k\geq\mu(\ba)=m+1$ for all $k$, so each $b_k=a_k-m\geq1$, which gives $\bb\in\P^K$.  Again by definition, we have $a_k=\mu(\ba)=m+1$ for some $k$, so for this $k$ we have $b_k=a_k-m=1$, and it follows from Lemma \ref{l1} that $\bb$ is an atom.

For uniqueness, suppose $\ba=m\bone+\bb=n\bone+\bc$ for $m,n\in\N$ and $\bb,\bc\in\A$.  We must show that $m=n$ and $\bb=\bc$.  Without loss of generality, we assume that $m\leq n$.  Now, $\bb=(n-m)\bone+\bc$ with $n-m\geq0$, so since $\bb$ is an atom it follows that $n-m=0$, i.e.~$m=n$.  But then also $\bb=0\bone+\bc=\bc$.
\epf

The atoms of any semigroup are of course contained in any generating set for the semigroup.  Conversely, it follows directly from Lemma \ref{l2} that $\P^K$ is generated by its atoms, so we have proved the following:

\begin{prop}\label{p1}
The set $\A$ of all atoms is the (unique) minimum generating set for $\P^K$.  \qed
\end{prop}

It is therefore very natural to look for a presentation for $\P^K$ in terms of the generating set~$\A$.  To avoid notational clashes, we define an abstract alphabet
\[
X = \set{x_\ba}{\ba\in\A},
\]
in one-one correspondence with $\A$.  By Proposition \ref{p1}, we have a surmorphism
\[
\phi:X^+ \to \P^K:x_\ba\mt\ba.
\]
Let $R\sub X^+\times X^+$ denote the set of all relations of the form
\begin{equation}\label{R}
x_\ba x_\bb = x_\bone^m x_\bc \qquad\text{for $\ba,\bb\in\A$, where $m=\mu(\ba+\bb)-1$ and $\bc=\ba+\bb-m\bone$.}
\end{equation}

\begin{rem}\label{r1}
As in the proof of Lemma \ref{l2}, we have $m\in\N$ (in fact, $m\in\P$) and $\bc\in\A$ in~\eqref{R}, so it follows that $x_\bone^m x_\bc$ is indeed a word over $X$.

A couple of special cases of \eqref{R} are worth noting.  If $\ba=\bone$, then since $\bb$ is an atom, we have $m=1$ and $\bc=\bb$, so \eqref{R} merely says $x_\bone x_\bb=x_\bone x_\bb$ in this case.  When $\bb=\bone$,~\eqref{R} says $x_\ba x_\bone=x_\bone x_\ba$.
\end{rem}

Our goal now is to show that $\P^K$ has presentation $\pres XR$ via $\phi$.  For the rest of the paper, we write ${\sim}=R^\sharp$ for the congruence on $X^+$ generated by $R$.

\begin{lemma}\label{l3}
For any word $w\in X^+$ we have $w\sim x_\bone^mx_\ba$ for some $m\in\N$ and $\ba\in\A$.
\end{lemma}

\pf
We prove the lemma by induction on the length $k$ of $w$.  If $k=1$, then $w=x_\ba$ for some~$\ba\in\A$, and we are done (with $m=0$).  Now suppose $k\geq2$, and write $w=x_{\ba_1}\cdots x_{\ba_k}$ where $\ba_1,\ldots,\ba_k\in\A$.  By induction we have $x_{\ba_1}\cdots x_{\ba_{k-1}}\sim x_\bone^nx_\bb$ for some $n\in\N$ and $\bb\in\A$.  Since $\bb$ and $\ba_k$ are both atoms, $R$ contains a relation of the form $x_\bb x_{\ba_k} = x_\bone^qx_\ba$, for some $q\in\N$ and $\ba\in\A$.  But then
\[
w = x_{\ba_1}\cdots x_{\ba_{k-1}}x_{\ba_k} \sim x_\bone^n x_\bb x_{\ba_k} \sim x_\bone^n x_\bone^qx_\ba = x_\bone^mx_\ba,
\]
where $m=n+q\in\N$.
\epf

We are now ready to prove our main result.

\begin{thm}\label{t1}
The semigroup $\P^K$ has presentation $\pres XR$ via $\phi$.
\end{thm}

\pf
We have already noted that $\phi$ is a surmorphism, so it remains to show that $\ker\phi=R^\sharp$.  

To show that $R^\sharp\sub\ker\phi$, we must show that $\phi$ preserves each relation from $R$.  To do so, let $\ba,\bb\in\A$, and let $m$ and $\bc$ be as in \eqref{R}.  Then
\[
\phi(x_\bone^mx_\bc) = m\bone+\bc = m\bone+(\ba+\bb-m\bone) = \ba+\bb = \phi(x_\ba x_\bb).
\]

For the reverse inclusion, let $(u,v)\in\ker\phi$, meaning that $u,v\in X^+$ and $\phi(u)=\phi(v)$.  We must show that $u\sim v$.  By Lemma \ref{l3} we have $u\sim x_\bone^mx_\ba$ and $v\sim x_\bone^nx_\bb$ for some $m,n\in\N$ and $\ba,\bb\in\A$.  But then (using ${\sim}=R^\sharp\sub\ker\phi$)
\[
m\bone+\ba = \phi(x_\bone^mx_\ba) = \phi(u) = \phi(v) = \phi(x_\bone^nx_\bb) = n\bone+\bb,
\]
and it follows from Lemma \ref{l2} that $m=n$ and $\ba=\bb$.  Putting everything together, we now have
\[
u \sim x_\bone^mx_\ba = x_\bone^nx_\bb \sim v,
\]
and the proof is complete.
\epf

\begin{rem}\label{r2}
To conclude, we consider the special case of $\P^2=\set{(a,b)}{a,b\in\P}$.  Here we have
\[
\A = \set{(a,1)}{a\in\P} \cup \set{(1,b)}{b\in\P}.
\]
To simplify notation, we relabel the letters from $X$:
\[
x = x_{(1,1)} \COMMA  y_a = x_{(a,1)} \AND z_a = x_{(1,a)} \qquad\text{for $a\geq2$.}
\]
So now $X=\{x\}\cup\set{y_a,z_a}{a\geq2}$.
Keeping Remark \ref{r1} in mind, the relations from $R$ become the following, with $a,b\geq2$ in each relation:
\begin{gather*}
xy_a = y_ax \COMMA xz_a = z_ax \COMMA
y_ay_b = xy_{a+b-1} \COMMA z_az_b = xz_{a+b-1} , \\[2mm]
y_az_b = z_by_a = \begin{cases}
x^az_{b-a+1} &\text{if $a<b$}\\
x^{a+1} &\text{if $a=b$}\\
x^by_{a-b+1} &\text{if $a>b$.}
\end{cases}
\end{gather*}

(Theorem \ref{t1} also holds when $|K|=1$, where $\P^K=\P$.  Here $\A=\{1\}$, and denoting the unique element of $X$ by $x$, the only relation in $R$ is the trivial $xx=xx$, so we obtain the usual presentation~$\pres x\es$ for $\P$.)
\end{rem}

\footnotesize
\def\bibspacing{-1.1pt}
\bibliography{biblio}

\begin{thebibliography}{10}

\bibitem{Artin1947}
E.~Artin.
\newblock Theory of braids.
\newblock {\em Ann. of Math. (2)}, 48:101--126, 1947.

\bibitem{CR2020}
A.~Clayton and N.~Ru\v{s}kuc.
\newblock On the number of subsemigroups of direct products involving the free
  monogenic semigroup.
\newblock {\em J. Aust. Math. Soc.}, 109(1):24--35, 2020.

\bibitem{CPbook}
A.~H. Clifford and G.~B. Preston.
\newblock {\em The algebraic theory of semigroups. {V}ol. {I}}.
\newblock Mathematical Surveys, No. 7. American Mathematical Society,
  Providence, R.I., 1961.

\bibitem{East2020}
J.~East.
\newblock Presentations for tensor categories.
\newblock {\em Preprint}, 2020, {\tt arXiv:2005.01953}.

\bibitem{Howie}
J.~M. Howie.
\newblock {\em Fundamentals of semigroup theory}, volume~12 of {\em London
  Mathematical Society Monographs. New Series}.
\newblock The Clarendon Press, Oxford University Press, New York, 1995.
\newblock Oxford Science Publications.

\bibitem{HR1994}
J.~M. Howie and N.~Ru\v{s}kuc.
\newblock Constructions and presentations for monoids.
\newblock {\em Comm. Algebra}, 22(15):6209--6224, 1994.

\bibitem{J1980}
D.~L. Johnson.
\newblock {\em Topics in the theory of group presentations}, volume~42 of {\em
  London Mathematical Society Lecture Note Series}.
\newblock Cambridge University Press, Cambridge-New York, 1980.

\bibitem{Jones1983_2}
V.~F.~R. Jones.
\newblock Index for subfactors.
\newblock {\em Invent. Math.}, 72(1):1--25, 1983.

\bibitem{L1998}
T.~G. Lavers.
\newblock Presentations of general products of monoids.
\newblock {\em J. Algebra}, 204(2):733--741, 1998.

\bibitem{LZ2015}
G.~Lehrer and R.~B. Zhang.
\newblock The {B}rauer category and invariant theory.
\newblock {\em J. Eur. Math. Soc. (JEMS)}, 17(9):2311--2351, 2015.

\bibitem{MKS1966}
W.~Magnus, A.~Karrass, and D.~Solitar.
\newblock {\em Combinatorial group theory: {P}resentations of groups in terms
  of generators and relations}.
\newblock Interscience Publishers [John Wiley \& Sons, Inc.], New
  York-London-Sydney, 1966.

\bibitem{MR2018}
P.~Mayr and N.~Ru\v{s}kuc.
\newblock Finiteness properties of direct products of algebraic structures.
\newblock {\em J. Algebra}, 494:167--187, 2018.

\bibitem{MR2019}
P.~Mayr and N.~Ru\v{s}kuc.
\newblock Generating subdirect products.
\newblock {\em J. Lond. Math. Soc. (2)}, 100(2):404--424, 2019.

\bibitem{RRW1998}
E.~F. Robertson, N.~Ru\v{s}kuc, and J.~Wiegold.
\newblock Generators and relations of direct products of semigroups.
\newblock {\em Trans. Amer. Math. Soc.}, 350(7):2665--2685, 1998.

\end{thebibliography}
\bibliographystyle{abbrv}

\end{document}